\makeatletter
\def\@maketitle{%
\defaultfont\normalsize
\let\@makefnmark\relax \let\@thefnmark\relax \ifx\@empty\@subjclass\else
\@footnotetext{1991 {\it Mathematics Subject
Classification}.\enspace
\@subjclass.}\fi
\ifx\@empty\@keywords\else
\@footnotetext{{\it Key words and phrases.}\enspace \@keywords.}\fi
\ifx\@empty\@thanks\else
\@footnotetext{\@thanks}\fi
\topskip66\p@ % 6.5 picas to the base of the first title line
\vtop{\centering{\baselineskip14\p@\bf
\expandafter{\@title}\@@par}%
\global\dimen@i\prevdepth}%
\prevdepth\dimen@i
\ifx\@empty\@authors
\else
\baselineskip32\p@
\vtop{\@andify{ AND }\@authors
\centering{{\@authors}\@@par}%
\global\dimen@i\prevdepth}\relax
\prevdepth\dimen@i
\fi
\ifx\@empty\@dedicatory
\else
\baselineskip18\p@
\vtop{\centering{\small\it\@dedicatory\@@par}%
\global\dimen@i\prevdepth}\prevdepth\dimen@i \fi
\ifx\@empty\@date\else
\baselineskip24\p@
\vtop{\centering\@date\@@par
\global\dimen@i\prevdepth}\prevdepth\dimen@i \fi
\normalsize
\dimen@32\p@ \advance\dimen@-\baselineskip \vskip\dimen@\@plus14\p@
} % end \@maketitle
\makeatother

\documentstyle[amscd,amssymb,verbatim,amsthm,12pt]{amsart}

\theoremstyle{plain}
\newtheorem{Thm}{Theorem}
\newtheorem*{Thm*}{Theorem}
\newtheorem{Cor}[Thm]{Corollary}
\newtheorem{Lem}[Thm]{Lemma}
\newtheorem{Prop}[Thm]{Proposition}
\newtheorem{Def}[Thm]{Definition}

\theoremstyle{remark}
\newtheorem{Ex}[Thm]{Example}
\newtheorem{Rem}[Thm]{Remark}

\newcommand{\ra}{\rightarrow}
\newcommand{\Spec}{\operatorname{Spec}}

\renewcommand{\O}{{\Cal O}}

\def\uint#1#2{{U\kern -8pt\int}_{\kern -5pt {#1}}^{\kern 2pt {#2}} }
\def\lint#1#2{{L\kern -7.5pt \int}_{\kern -5pt {#1}}^{\kern 2pt {#2}} }
\def\morph#1{\overset{#1}{\ra}}

\def\dsp#1{$\displaystyle{#1}$}
\def\Cal#1{{\cal #1}}
\def\BBB#1{{\Bbb #1}}

\theoremstyle{remark}

\numberwithin{equation}{section}

\newcommand{\Prod}{\operatornamewithlimits{\Pi}}

\newcommand{\A}{{\Cal A}}

\renewcommand{\O}{{\Cal O}}
\newcommand{\F}{{\Cal F}}

\newcommand{\D}{{\Cal D}}

\newcommand{\Pic}{\operatorname{Pic}}
\newcommand{\ann}{\operatorname{ann}}

\newcommand{\Proj}{\operatorname{Proj}}

\newcommand{\C}{\BBB C}
\newcommand{\Z}{\BBB Z}

\def\Mod{\text{-mod}}

\begin{document}

\title{$\D$-modules on $1|1$ supercurves}

\author{Mitchell J. Rothstein} 
\address{Department of Mathematics \\ University of Georgia \\ }
\email{rothstei@@math.uga.edu}

\author{Jeffrey M. Rabin}
\address{Department of Mathematics \\ UCSD \\ La Jolla, CA 92093}
\email{jrabin@@ucsd.edu}

\maketitle
%*********************************************************************

\begin{abstract}
It is known that to every $(1|1)$ dimensional supercurve  $X$ there is associated a dual supercurve $\hat X$,  and a superdiagonal  $\Delta\subset X\times \hat X$.  We establish that the categories of $\D$-modules on $X$,  $\hat X$ and $\Delta$ are equivalent.  This follows from a more general result about $\D$-modules and purely odd submersions.  The equivalences  preserve tensor products,  and take vector bundles to vector bundles.  Line bundles with connection are studied,  and examples are given where $X$ is a superelliptic curve.
\end{abstract}

\section{Introduction}

Super Riemann surfaces, the analogue of Riemann surfaces in the category of supermanifolds, have been extensively studied since their introduction in the context of superstring theory \cite{CR1988, F1986, LR1988, RSV1988}.
They are supercurves of dimension $(1|1)$ satisfying an additional ``superconformal" constraint which allows one to identify irreducible (Weil) divisors with points.  The study of more 
general supercurves $X$, without the constraint,  began in earnest with the  paper \cite{DRS1990},  where it was observed that  the irreducible divisors on $X$ can be identified with points of a dual supercurve, $\hat X$, having the same underlying topological space, and that the dual of $\hat X$ is again $X$.   Furthermore,  there is a distinguished $(1|2)$-dimensional submanifold,
$$\Delta\subset X\times\hat X$$
(the ``superdiagonal"), which exhibits $\hat X$ as the family of irreducible divisors on $X$ and vice versa.    Supercurves,  and their duality,  were later seen in \cite{BR1999} to play an important role in the study of super analogues of the KP-hierarchy.  In that paper the main objects of  study were $\Pic(X)$ and $\Pic(\hat X)$.

In this paper we explore the categories  of $\D$-modules over $X$, $\hat X$ and $\Delta$.   We find in section 4 that in fact the three categories $\D_X$-mod,  $\D_{\hat X}$-mod and $\D_{\Delta}$-mod are equivalent.
This follows from a more general statement  in section 2,  about $\Cal D$-modules and purely odd submersions.  Section 3 provides a review of supercurve duality.   The equivalence of categories is further explored in sections 5--7, where explicit formulas in local coordinates are given.  Section 8 describes the direct image of a trivial vector bundle with connection, and section 9 specializes this to line bundles. Section 10 illustrates our results with explicit examples for the case of super elliptic curves, i.e., supercurves of genus one.

Throughout the paper we work in the category of supermanifolds over a fixed superscheme $S$ over $\BBB C$.  By supermanifold we mean a smooth morphism $Z\to S$, and by dimension we mean the relative dimension. Though most of the results in this paper are valid for arbitrary $S$,  we will avoid notational complications by assuming that $S=\Spec(\Lambda)$,  where $\Lambda$ is a finite-dimensional nilpotent extension of $\BBB C$.  Except as noted,  the results are valid in both the Zariski and complex topologies.

\section{$\Cal D$-modules and purely odd submersions}
Let $\sigma:Z\to W$ be a  smooth morphism of smooth superschemes over $S$.  
Let $\D_Z$  denote the sheaf of linear differential operators on $\O_Z$.  Then a $\D_Z$-module is a sheaf $\F$ of $\O_Z$-modules  equipped with a flat connection $$\nabla:\Omega_Z\otimes \F\to\Omega_Z[1]\otimes\F$$
where $\Omega_Z$ is the sheaf of differential one-forms relative to $S$, and $\Omega_Z[1]$ denotes $\Omega_Z$ with a degree shift.
(For background on $\D$-modules,  see \cite{Bj1993,B1987}.)
One has direct and inverse image functors for $\D$-modules.  The inverse image of a $\D_W$-module  $\F$ has $\sigma^*(\F)$ as its underlying $\O_Z$-module.  The direct image is, in general,  defined in the derived category.   

Assume now that $\sigma$ is a submersion.   Then there is an underived version of the direct image, defined as follows. Let  $\Cal T_{Z}$ denote the  tangent sheaf  of $Z$ and let $\Cal T_{\sigma}\subset\Cal T_Z$ denote the vertical tangent sheaf.  Then we have an exact sequence
$$0\to\Cal T_{\sigma}\to\Cal T_Z\to \sigma^*(\Cal T_W)\to 0$$
The direct image functor
$$\sigma_+:\D_Z\Mod\to\D_W\Mod$$
is defined  by
\begin{equation}\label{eq:pushforward}\sigma_+(\F)= \sigma_*(\ann(\Cal T_{\sigma}) )\end{equation}
where $\ann$ denotes annihilator.  

Say that the submersion $\sigma:Z\to W$ is {\em purely odd} if the fibers have dimension 
$(0|n)$ for some $n$.  The key observation in the paper is the following result.

\begin{Thm}\label{thm:oddEquivalent}
Let $\sigma:Z\to W$ be a purely odd submersion.  
Then the categories of $\D$-modules on $Z$ and $W$ are equivalent.  Specifically,  the functors 
$\sigma^*$ and $\sigma_+$ are inverse equivalences.
\end{Thm}

\begin{pf} To simplify the notation,  note that $Z$ and $W$ share the same underlying topological space.
If $\Cal G$ is a $\D_W$-module,  then $\Cal G$ maps naturally to $\sigma^*(\Cal G)$. 
A computation in local coordinates easily shows that this map is injective,  and that
the image of $\Cal G$ is precisely $\sigma_+\sigma^*(\Cal G)$.  

If $\F$ is a $\D_Z$-module,
then $\sigma_+({\Cal F})\subset{\Cal F}$,  and we have a natural map
$$\O_Z\otimes_{\O_W}\sigma_+({\Cal F})=\sigma^*\sigma_+({\Cal F})\to {\Cal F}$$
The fact that this map is an isomorphism follows from the purely algebraic lemma stated below.\end{pf}

%Theorem \ref{thm:oddEquivalent} is a corollary of the following purely algebraic lemma:

\begin{Lem}\label{lem:anticomm}
Let $R=R_0\oplus R_1$ be a $\BBB Z_2$-graded ring.   Let
$$Q=R[\theta_1,...,\theta_n,\partial_1,...,\partial_n]\ $$
where the $\theta_i$'s are free supercommuting odd variables, $\partial_i=\partial/\partial\theta_i$, and the
$\theta$'s and $\partial$'s supercommute with $R$. Then the categories $Q$-mod and $R$-mod are equivalent.  Specifically,  the following functors are inverses:
\begin{align}
Q\Mod\ni M&\mapsto M_*=\ann(\partial_1,...,\partial_n)\\
R\Mod\ni N&\mapsto N^*= R[\theta_1,...,\theta_n]\otimes_RN
\end{align}
The statement that the natural homomorphism
$$(M_*)^*\simeq M
$$
is an isomorphism is equivalent to the statement that every $A\in M$ has a unique expansion
\begin{equation}\label{eq:expansion}
A=\sum_{\mu}A_{\mu}\theta^{\mu}\ , \  A_{\mu}\in M_*
\end{equation}
where $\mu=(\mu_1,\dots,\mu_n)$ is a multiindex of $0$'s and $1$'s.\end{Lem}

\begin{pf}
It is easy to see that for all $R$-modules $N$,  $N\simeq (N^*)_*$.   

For the other direction,  
one has the natural map
\begin{equation}R[\theta_1,\dots,\theta_n ]\otimes_R M_*\to M
\end{equation}
which is an isomorphism if and only if formula \eqref{eq:expansion} is valid.  We prove \eqref{eq:expansion} by induction on $n$.  For $n=1$,   
let $A\in M$.   Then the decomposition is given by
 \begin{align}\label{eq:decomp}
A_0&=A-\theta_1\partial_1A\\A_1&=\partial_1A
\end{align}
For the uniqueness,  apply $\partial_1$ to both sides of the equation

 $$A=A_0+\theta_1A_1$$

If $n>1$,  write
$$R[\theta_1,...,\theta_n,\partial_1,...,\partial_n]=R[\theta_2,...,\theta_n,\partial_2,...,\partial_n][\theta_1,\partial_1]$$

\end{pf}

\begin{Prop}\label{prop:preserveRank}
Let $\sigma:Z\to W$ be a purely odd submersion.  Then $\sigma_+$  is exact, and preserves tensor products.  It
takes $\D_Z$-modules that are locally free as $\O_Z$-modules to locally free  $\O_W$-modules, and  preserves rank.    
\end{Prop}

\begin{pf} The first part follows immediately from theorem \ref{thm:oddEquivalent}.    Since the proposition is local,  we may assume for the second part that we have a flat connection on the trivial bundle of rank $p|q$,  i.e., $\O_Z^{p|q}$.  Let the fiber dimension of $\sigma$ be $0|n$,  and let $\theta_1,...,\theta_n$ be a set of fiber coordinates.
%
% A connection on $\O_Z^{p|q}$ is of the form $d+\omega$, where $\omega$ is a  one-form with values in $\frak {gl}^{p|q}(\O_Z)$, satisfying the zero-curvature condition
%\begin{equation}\label{eq:flat}
%d\omega+\omega\wedge\omega=0.
%\end{equation}
Let $I\in \frak {gl}^{p|q}(\O_Z)$ denote the identity matrix.   Then we have a unique decomposition
$$I=\sum_{\mu}A_{\mu}\theta^{\mu}$$
where, for all $i$,  $\nabla_{{\theta_i}} A_{\mu}=0$.  Let $A=A_{(0,\dots,0)}$. Then $A$  is an invertible matrix. Thus the columns of $A$ lie in $\sigma_+(\O_Z^{p|q})$ and form a basis for $\O_Z^{p|q}$. If $\psi\in\sigma_+(\O_Z^{p|q})$, then there is a unique  vector
$\phi\in \O_Z^{p|q}$ such that $\psi=A\phi$. Then $0=\nabla_{\theta_i}(A \phi)=A\partial_{\theta_i}\phi$, whence the entries of $\phi$ belong to $\O_W$.  Thus the columns of $A$ form a basis for $\sigma_+(\O_Z^{p|q})$ as an $\O_W$-module.
\end{pf}

\begin{Rem}The standard result in the commutative setting, that a $\D$-module is locally free of finite rank as an $\O$-module if and only if it is coherent as an $\O$-module \cite{Bj1993,B1987}, holds in the supercommutative setting as well.\end{Rem}

\section{Supercurves and their duals}

 By definition,  a  supercurve is a supermanifold of dimension $(1|n)$ for some $n\ge 1$. Let $X$ be a supercurve of dimension $(1|1)$.
 Then there is a dual  $(1|1)$-dimensional supercurve  $\hat X$ constructed as follows \cite{DRS1990,BR1999}.   Define
\begin{equation}\label{eq:exterior}
\Delta_X=\Proj(\Omega_X)\morph{\pi} X\end{equation}
Thus,  if $(z,\theta)$ are local coordinates on an open set $\Cal U$, $\Omega_X(\Cal U)$ is the polynomial algebra $\O_X(\Cal U)[d\theta,dz]$,  with $d\theta$ even and $dz$ odd.  In the proj construction,  $d\theta, dz$ are taken as homogeneous coordinates,  of which only $d\theta$ may be inverted.  Thus,  on $\pi^{-1}(\Cal U)$ one has the local coordinate system
\begin{equation}\label{eq:coord} ( z,\theta,\rho) ,\end{equation}
where
\begin{equation}\label{eq:rho}\rho=d\theta^{-1}dz\ .\end{equation}
In particular,  $\Delta_X$ has relative dimension $(1|2)$ over $S$.
%(The reason for the somewhat pedantic underlining of the variables is that $dz$ and $d\theta$ are sections of $\O_{\Delta_X}(1)$,  and therefore need to be distinguished from $d z$ and $d\theta$, which are sections of $\Omega_{\Delta_X}$.)

The exterior derivative $d$  is an odd derivation of degree one: $$d:\Omega_X\to\Omega_X[1]\ .$$
By the standard sheaf construction \cite{H1977},  one obtains an odd derivation
\begin{equation}\label{eq:dtilde}
\tilde d:\O_{\Delta_X}\to\O_{\Delta_X}(1)\ .
\end{equation}
Note that $d\theta$ is a trivialization of $\O_{\Delta_X}(1)$.  Then we have the formula
 
\begin{equation}\label{eq:dtildefmla}
\tilde d=dz\ \partial_{ z}+d\theta\ \partial_{\theta}=d\theta(\rho\partial_{ z}+\partial_{ \theta})\ .\end{equation}
The {\em dual curve},  $\hat X$, has the same topological space as $X$, with structure sheaf
\begin{equation}\label{eq:Xhat}
\O_{\hat X}=\ker(\tilde d)\  .\end{equation}
Let 
$$ u= z-\theta\rho\ .$$
%Adopting the same naming convention for $\O_{\hat X}$,  
Then $\tilde d\rho=\tilde d u=0$.  Furthermore,  $( u,\rho,\theta)$ is a local coordinate system on $\Delta_X$,  and one checks that 
  $(u,\rho)$ is a local coordinate system on $\hat X$.  In particular,  $\hat X$ is a family of smooth $(1|1)$-supercurves over  $S$.  
It is known \cite{BR1999} that $\Delta_{ X}$ and $\Delta_{\hat X}$ are naturally isomorphic, as superschemes over $\hat X$.  On the level of structure sheaves,  the isomorphism is given in local coordinates by
   \begin{align}\label{eq:Delta}
  & \O_{\Delta_{\hat X}}\to\O_{\Delta_X}\notag\\
  &u\mapsto z-\theta\rho\\
  &\rho\mapsto\rho\\
  &d\rho^{-1}{du}\mapsto \theta
  \end{align}
  Then $\O_{\hat{\hat X}}$ appears as a subsheaf of $\O_{\Delta_X}$,  and one checks that this subsheaf coincides with the image of $\O_X$.  Thus,
$\hat{\hat X}$ is naturally isomorphic to $X$.  In local coordinates, the isomorphism is given by
  \begin{align}
  &u-\rho\frac{du}{d\rho}\mapsto z\\
&d\rho^{-1}{du}\mapsto \theta
  \end{align}
  One should view $\Delta_{ X} \cong \Delta_{\hat X}$ (denoted simply $\Delta$ in the sequel) as the 
  ``superdiagonal" in $X \times {\hat X}$ defined by the equation $z-u-\theta\rho=0$.
  
  \section{Equivalences of categories}
  
 \begin{Thm}\label{thm:main}The categories $\D_X$-mod,  $\D_{\hat X}$-mod and 
 $\D_{\Delta}$-mod are equivalent.\end{Thm}
 
 \begin{pf}The maps $\pi:\Delta\to X$ and $\hat\pi:\Delta \to \hat X$ are purely odd submersions,  so the result follows from theorem \ref{thm:oddEquivalent}.\end{pf}
 
 For $\F$ a $\D_X$-module, define
 \begin{equation}\label{eq:hat}
\hat\F=\hat\pi_+\pi^*(\F)
\end{equation}

Then we have a canonical isomorphism  $\hat{\hat\F}\simeq F$,  by theorem \ref{thm:main}.

\begin{Ex}
By definition,  $\hat{\O}_X=\O_{\hat X}$.  

\end{Ex}
\begin{Ex}

Consider $\D_X$ as a left $\D_X$-module. We have
$$\pi^*(\D_X)=\text{Diff}(\O_X,\O_\Delta)\ .$$
Then a germ $L\in \text{Diff}(\O_X,\O_\Delta)$ belongs to $\hat{\D}_X$ if and only if
$${\tilde d}\circ L=0$$
where $\tilde d$ is as in  \eqref{eq:dtilde}.
We therefore have
$$\hat{\D}_X=\text{Diff}(\O_X,\O_{\hat X})$$
which is to say, differential operators from $\O_X$ to $\O_\Delta$ that factor through the inclusion
$\O_{\hat X}\to\O_\Delta$.

\end{Ex}

\section{Local description}

Let $\Cal U\subset X$ be an open set,  with coordinates $(z,\theta)$.  Let $(u,\rho)$ be the corresponding coordinates on $\hat X$.   Then we get an isomorphism
\begin{equation}\label{eq:localIso}
\O_X|_{\Cal U}\overset{\Psi^{(z,\theta)}}\longrightarrow\O_{\hat X}|_{\Cal U}
\end{equation}
sending $z\to u $ and $\theta\to\rho$.   That is,
\begin{align}\label{eq:IsoLocal}
\Psi^{(z,\theta)}(f(z)+\theta g(z))&=f(z-\theta\rho)+\rho g(z-\theta\rho)\\&=f(z)+\rho(\theta\partial_zf+g)\label{eq:notice}
\end{align}
The isomorphism $\Psi^{(z,\theta)}$ extends to an isomorphism
$$\Psi^{(z,\theta)}:\D_X|_{\Cal U}\overset{\sim}\to\D_{\hat X}|_{\Cal U}$$
 sending $\partial_z\mapsto \partial_u $ and $\partial_\theta\mapsto\partial_\rho$.  
This identification allows one to regard every  $\D_X|_{\Cal U}$-module as a $\D_{\hat X}|_{\Cal U}$-module  (in a coordinate-dependent way).  If ${\Cal F}$ is a  $\D_X|_{\Cal U}$-module,  let
 ${\Cal F}^{(z,\theta)}$ denote ${\Cal F}$ itself, regarded as a  $\D_{\hat X}|_{\Cal U}$-module.
 Notice that  equation \eqref{eq:notice}  can be written
 \begin{equation}\label{eq:written}
\Psi^{(z,\theta)}(h)=(1-\theta\partial_\theta+\rho(\theta\partial_z+\partial_\theta))(h)
\end{equation}
where $h=f+\theta g$.  This suggests  the following definition.
\begin{Def}\label{def:tau}Let
\begin{equation}\label{eq:tau}
\tau^{(z,\theta)}=1-\theta\nabla_\theta+\rho(\theta\nabla_z+\nabla_\theta)
\end{equation}
(Here it is understood that one has a $\D_X|_{\Cal U}$-module ${\Cal F}$,  and $\tau^{(z,\theta)}$ acts on $\pi^*{\Cal F}$.)
\end{Def}

\begin{Thm}\label{thm:localPicture}
Let  ${\Cal F}$ be  a $\D_X|_{\Cal U}$-module.  Then $\tau^{(z,\theta)}$  gives an isomorphism
\begin{equation}\label{eq:localPicture}
\tau^{(z,\theta)}:{\Cal F}^{(z,\theta)}\overset\sim\to\hat{{\Cal F}}
\end{equation}
\end{Thm}

\begin{pf}
Let us establish first that  $\tau^{(z,\theta)}$ gives a bijection ${\Cal F}\to\hat{{\Cal F}}$.
In $(z,\theta,\rho)$ coordinates,  $\Cal T_{\hat\pi}$ is generated by
     \begin{equation}\label{eq:boldD}{\bold D}= \rho\partial_z+\partial_{\theta}\end{equation}
   Sections of $\pi^*\F$ may be uniquely written in the form  
    $\phi=P+\rho Q$,  where $P$ and $Q$ are  sections of $\F$.
  Then  $\phi\in\hat \F$ if and only if $\nabla_{\bold D} \phi=0$, that is,
  \begin{equation}\label{eq:condition}
  \nabla_{\theta}P+\rho(\nabla_z P-\nabla_{\theta}Q)=0\ .\end{equation}
Set $P=P_0+\theta P_1$ and $Q=Q_0+\theta Q_1$,  where $\nabla_{\theta}$ annihilates $P_i$ and $Q_i$.
Then equation \eqref{eq:condition} reads
\begin{equation}
  P_1+\rho(\nabla_zP_0+\theta\nabla_zP_1-Q_1)=0\ .\end{equation}
  That is,
   \begin{equation}
  P_1=0\  ,\  \nabla_zP_0=Q_1\ .\end{equation}
  Then
   \begin{equation}
  \phi=P_0+\rho(Q_0+\theta\nabla_zP_0)= \tau^{(z,\theta)}(P_0+\theta Q_0)\ .\end{equation}
Since $\phi=0$ if and only if $P_0=Q_0=0$, $\tau^{(z,\theta)}$ is bijective.

It remains to prove that for all  $M\in\D_X|_{\Cal U}$ and all $A\in{\Cal F}$,
\begin{equation}\label{eq:toProve}
\tau^{(z,\theta)}(MA)=\Psi^{(z,\theta)}(M)\tau^{(z,\theta)}(A)
\end{equation}
We leave it to the reader to check this when $M\in\Cal O_X$.  It remains to check equation \eqref{eq:toProve} when $M$ is a partial derivative.
 To distinguish between partial derivatives in the $(z,\theta,\rho)$ coordinate system and the the $(u,\rho,\theta)$ coordinate system,  we will denote the latter by $\hat\partial_u$,  etc..  We also write $\hat\nabla_{u}$, etc..
 We then have
 \begin{align}
\hat\nabla_u&=\nabla_z\label{eq:u}\\
\hat\nabla_{\rho}&=\nabla_{\rho}-\theta\nabla_z\label{eq:rho}\\
\hat\nabla_{\theta}&=\nabla_{\theta}-\rho\nabla_z\label{eq:theta}
 \end{align}
 (These operators  are acting on $\pi^*\F$.)
% 
% Notice that
% \begin{equation}\label{eq:tau}
% \tau^{z,\theta}=1-\theta\nabla_{\theta}+\rho(\nabla_{\theta}+\theta\nabla_z)\end{equation}
 Then $$\hat\nabla_u\tau^{(z,\theta)}=\nabla_z\tau^{(z,\theta)}=\tau^{(z,\theta)}\nabla_z$$ 
 We also have
 $$\nabla_{\rho}\tau^{(z,\theta)}=\tau^{(z,\theta)}\nabla_{\rho}+\nabla_{\theta}+\theta\nabla_z$$
 Then
 $$\hat\nabla_{\rho}\tau^{(z,\theta)}=\tau^{(z,\theta)}\nabla_{\rho}+\nabla_{\theta}+\rho\theta\nabla_z(\nabla_{\theta}+\theta\nabla_z)=\tau^{(z,\theta)}\nabla_{\rho}+\tau^{(z,\theta)}\nabla_{\theta}$$
 The result follows,  since one is applying this operator to the kernel of $\nabla_{\rho}$.
   \end{pf}
   
  Theorem \ref{thm:localPicture} provides a \v Cech description of  the functor ${\Cal F}\to \hat{\Cal F}$.  Cover $X$ by coordinate charts $(\Cal  U_i,z_i,\theta_i)$.  Then $\O_{\hat X}$ may regarded as $\O_X$ glued to itself by transition automorphisms $D_{i,j}={\Psi^{(z_i,\theta_i)}}^{-1}\Psi^{(z_j,\theta_j)}$.  Then $D_{i,j}$ is a differential operator,  so it can be applied to an arbitrary $\D_X$-module ${\Cal F}$.  Then theorem \ref{thm:localPicture} shows that $\hat {\Cal F}$ is ${\Cal F}$ glued to itself by the cocycle $D_{i,j}$.

\section{Projected and injected supercurves}
A {\em projected supercurve over $S$,}   \cite{BR1999},  is a submersion  $\sigma:X\to X_0$ of smooth superschemes over $S$,  where $X$ and $X_0$ have relative dimension $(1|1)$ and $(1|0)$ over $S$, respectively.    With the same assumptions on $X$ and $X_0$,  an {\em injected} supercurve over $S$ is an immersion $\iota:X_0\to X$.

Let $\sigma:X\to X_0$ be a projected supercurve. 
%Now let $X_0/S$ be a smooth curve over $S$.   That is,  the relative dimension of $X_0$ over $S$ is $(1|0)$.
%If $X$ is our supercurve,  and if there is a submersion $\sigma:X\to X_0$,  then 
By theorem
\ref{thm:oddEquivalent},  the categories $\D_X\Mod$ and $\D_{X_0}\Mod$ are equivalent.
%supercurves equipped with such a submersion are called {\em projected supercurves} in \cite{BR1999}.
The following result is proved in \cite{OR2002}.  For the reader's  convenience we give a proof here.

\begin{Prop}\label{prop:switcharoo}
Fix a smooth curve $X_0/S$.  Then the category of  projected  supercurves
$\sigma:X\to X_0$ is equivalent to the category of  injected
supercurves
$\iota:X_0\to X$.  The equivalence is given by
$X\mapsto \hat X$. \end{Prop}

 \begin{pf} Let $\sigma:X\to X_0$ be a projected supercurve.  Then we have $\O_{X_0}\subset\O_X$.
% 
% Then we have
% \begin{align}\sigma^*:{\Omega_{X_0}}_{cl}&\to {\Omega^1_X}=\O_{\Delta X}(1)
%\end{align}
%where the $cl$ subscript means closed.   Thus
Let $(z,\theta)$ and $(w,\eta)$ be two local coordinate systems on an  open subset, such that
$z,w\in \O_{X_0}$.  Let $\rho=d\theta^{-1}dz$,  $\lambda=d\eta^{-1}dw$.  Writing
$w=f(z)$ and $\eta=\theta g(z)+\Lambda(z)$,  we have
\begin{align}
\lambda&=\frac{\rho f'(z)}{\rho (\theta g'(z)+\Lambda'(z))+g(z)}\\
&=\frac{\rho f'(z)}{g(z)} 
\end{align}
Thus we have a globally defined ideal $\Cal I\subset\O_\Delta$,  spanned locally by $\rho$. 
It is easily seen that we have an exact sequence
$$0\to\Cal I\to \O_\Delta \overset{\alpha}\to \O_X\to 0$$
Consider the restriction of $\alpha$ to $\O_{\hat X}$.   Write a section of $\O_\Delta$ as
$A+\rho B$,  $A, B\in\O_X$.    According to equation  \eqref{eq:condition},  $A+\rho B\in \O_{\hat X}$ if and only if
$\partial_{\theta} A=0$,  $\partial_z A=\partial_{\theta} B$.  In particular,  $A\in\O_{X_0}$.   Furthermore,  for all $A\in \O_{X_0}$,  we have $A+\rho\theta\partial_zA\in\O_{\hat X}$.  Thus,
$\alpha$ restricts to a surjection $\O_{\hat X}\to\O_{X_0}$, or equivalently, an immersion $X_0\to\hat X$.

Conversely,  let $\iota:X_0\to X$ be an injected supercurve.  Let $\Cal J$ denote the kernel of the corresponding surjection $\iota^*:\O_X\to\O_{X_0}$.  Let $\Cal U$ be an open set  sufficiently small that   $\Cal J|_{\Cal U}$ is generated by an odd section $\theta$. Then
$d\theta$ trivializes $\O_\Delta(1)|_{\Cal U}$.  We get a homomorphism of sheaves of rings,
\begin{align}
\Cal O_X|_{\Cal U} &\to \Cal O_{\hat X}|_{\Cal U}\\
f&\mapsto \frac{d(\theta f)}{d\theta}=f-\theta\frac{d f}{d\theta}\label{eq:nicemap}
\end{align}
(This is the piece of $\tau^{(z,\theta)}$ that doesn't depend on a choice of $z$.)
Let $\eta$ be  another generator of $\Cal J|_{\Cal U}$,  and let $\eta=g\theta$, where $g\in\O_X^*(\Cal U)$.  Then
\begin{align}
 \frac{d(\eta f)}{d\eta}&=f-g\theta\frac{d f}{gd\theta-\theta dg} \\
 &=f-g\theta\frac{d f}{gd\theta} =\frac{d(\theta f)}{d\theta}
\end{align}
Thus we have a globally defined homomorphism of sheaves of rings,
\begin{equation}\label{eq:homo}\nu:\Cal O_X \to \Cal O_{\hat X}\end{equation}
If we now extend $\theta$ to a coordinate system $(z,\theta)$,  then
$$\nu(f(z)+\theta g(z))=f+\rho\theta\partial_z f$$
Thus the kernel of $\nu$ is $\Cal J$, and we have produced an injection $$\O_{X_0}\to \O_{\hat X}$$
%Let $\Cal I'\subset$
%Write a section  $\phi\in \O_{\hat X}$ as
%\begin{align}
%\phi&=\tau^{(z,\theta)}(A+\theta B)\\
%&=\tau^{(w,\eta)}(\tilde A+\eta\tilde B)
%\end{align}
%where $A$, $B$, ..., belong to $\O_{X_0}$.
 \end{pf}
 
 With $X$ and $X_0$ as above,  let $\sigma:X\to X_0$ be a projected supercurve with corresponding injected supercurve  $\iota:X_0\to \hat X$.  We then have two functors $\D_X\Mod\to\D_{X_0}\Mod$,
 \begin{align}
 \F&\mapsto\iota^*(\hat\F)\label{eq:functor1}\\
 \F&\mapsto\sigma_+(\F)\label{eq:functor2}
 \end{align}

\begin{Prop}The functors (\ref{eq:functor1}) and (\ref{eq:functor2}) are naturally isomorphic.\end{Prop}
\begin{pf}  Let $\F$ be a $\D_X$-module.
Formula \eqref{eq:nicemap} can be used on $\F$,  to give a globally defined map
\begin{align}\label{eq:modulehomo}\nu:\F&\to \hat \F\\
\phi&\mapsto \frac{\nabla(\theta \phi)}{d\theta}\end{align}
Restrict $\nu$ to $\sigma_+(\F)$ and follow it with the pullback map $\hat \F\to\iota^*(\hat\F)$ to obtain the desired isomorphism.
\end{pf}

\section{Split Supercurves}  Continuing with $X$ and $X_0$ as above,  $X$ is said to be {\em split} over $X_0$ if there is a locally free rank-one sheaf ${\Cal L}$ of $\O_{X_0}$-modules such that the structure sheaf of $X$ is $\O_{X_0}\oplus{\Cal L}$.  This is equivalent to the existence of both a submersion
$\sigma:X\to X_0$ and an immersion $\iota:X_0\to X$,  such that $\sigma\circ\iota=id$.  Then by
proposition \ref{prop:switcharoo},  $\hat X$ is also split over $X_0$,  with the line bundle in question being the Serre dual of the original.   Then by theorem \ref{thm:oddEquivalent},   we can identify
the categories $\D_{X}$-mod and $\D_{\hat X}$-mod  with $\D_{X_0}$-mod.

\begin{Prop}\label{prop:trivial}Let $X$ be split over $X_0$. Under the equivalences
$$ \D_{X_0}\text{-mod}  \cong \D_X\text{-mod }\ ,\   \D_{X_0}\text{-mod}   \cong  \D_{\hat X}\text{-mod}$$
the transform $\F\to\hat \F$ reduces to the identity functor on $\D_{X_0}{ -mod}$.
\end{Prop}
\begin{pf}   Given the submersion $\sigma:X\to X_0$ and immersion $\iota:X_0\to X$,
it is easy to check that the two functors $\sigma_+$ and $\iota^*$ are naturally isomorphic.
Let $\Cal G$ be a $\D_{X_0}$-module and write $\Cal G=\sigma_+(\F)$. Then the functor we are considering sends $\Cal G$ to $\sigma_+(\hat\F)$.
We have
$$\sigma_+(\hat\F)\cong \iota^*(\F)\cong \sigma_+(\F)\cong\Cal G $$\end{pf}

\section{Direct image of the trivial bundle  with connection}
Returning now to the purely odd submersion $\sigma:Z\to W$,  consider a connection $d+\omega$ 
on the trivial bundle $\O_Z^{p|q}$.  Here
 $\omega$ is a  one-form with values in ${\frak {gl}}^{p|q}(\O_Z)$, satisfying the zero-curvature condition
\begin{equation}\label{eq:flat}
d\omega+\omega\wedge\omega=0
\end{equation}
According to proposition \ref{prop:preserveRank},  $\sigma_+(\O_Z^{p|q},d+\omega)$  is a locally free $\O_W$-module, of  rank $p|q$.
It is natural to ask for a description of  this $\O_W$-module.

The one-form $\omega$ restricts to a relative flat connection form on the fibers of $\sigma$.  Denote this restriction by $\omega_{\sigma}$.

Let $\Cal S^{p|q}$ denote the sheaf of flat connection forms on $Z$ and let $\Cal S_{\sigma}^{p|q}$ denote the sheaf of relative flat connection forms.  Let $d_{\sigma}$ denote the relative differential.
The sheaf
$Gl^{p|q}(\O_Z)$ maps to $\Cal S_{\sigma}^{p|q}$ by 
\begin{equation}\label{eq:connection}
A\mapsto  -d_{\sigma}A\cdot A^{-1}
\end{equation}

For a superscheme $Z$, let $\Cal N_Z\subset\O_Z$ denote the sheaf of nilpotents.
\begin{Prop}\label{prop:onto} 
The sheaf of subgroups $I+{\frak gl}^{p|q}(\Cal N_Z)$ maps surjectively to  $\Cal S_{\sigma}^{p|q}$.
\end{Prop}

\begin{pf}
Let $I$ denote the $(p+q)\times(p+q)$ identity matrix.  Let $\theta_1,...,\theta_n$ be fiber coordinates 
on an open set $\Cal U\subset Z$.  Decompose $I$ as in equation \eqref{eq:expansion}.
$$I=\sum_{\mu}\theta^{\mu}A_{\mu}$$
where $(d+\omega)_{\sigma}(A_{\mu})=0$.
In particular, 
\begin{equation}\label{eq:particular}
-(d_{\sigma}A_0)A_0^{-1}=\omega_{\sigma}
\end{equation}
  Furthermore,  the proof of  lemma \ref{lem:anticomm} shows that
$$A_0=\Prod_{i=1}^n(1-\theta_i\nabla_{\theta_i})(I)\in I+{\frak gl}^{p|q}(\Cal N_Z)$$
\end{pf}

If we regard  $\Cal S_{\sigma}^{p|q}$ as a sheaf of pointed sets,  where the $0$ one-form is the distinguished point, then the kernel of the map \eqref{eq:connection} is the sheaf
$Gl^{p|q}(\O_W)$.  We therefore have a connecting homomorphism
\begin{equation}\label{eq:connecting}
H^0(Z,\Cal S_{\sigma})\to H^1(W,  I+{\frak gl}^{p|q}(\Cal N_W))
\end{equation}
($Z$ and $W$ share the same topological space.)

\begin{Cor}\label{cor:connecting}
Let ${\Cal F}$ be a $\D_Z$-module with underlying $\O_Z$-module $\O^{p|q}$ and connection one-form
$\omega$.  Regarding $\sigma_+({\Cal F})$ simply as a vector bundle,   its class belongs to 
 $H^1(W,  I+{\frak gl}^{p|q}(\Cal N_W))$,  and that class is the image of $\omega_{\sigma}$ under the connecting homomorphism \eqref{eq:connecting}.
\end{Cor}

\begin{pf}
With the notation as in proposition \ref{prop:onto},  the columns of $A_0$ form a basis for $\sigma_+({\Cal F})|_{\Cal U}$.  The result then follows from equation  \eqref{eq:particular}.
\end{pf}

\section{Line bundles with connection} 
Here we apply the results of the previous section to the transform of the trivial bundle with connection on $ X$. Thus,  we choose an odd,  closed one-form $\omega$ on $X$, and consider the $\D_X$-module $\O_X^{\omega}=(\O_X,d+\omega)$. 
%Let ${\Cal L}$ be a $\D_X$-module,  locally free of rank one as an $\O_X$-module,  (or to be brief,  a line bundle with connection on $X$.)   Then $\hat{{\Cal L}}$ is a line bundle with connection on $\hat X$.  In this section we address the problem of understanding $\hat{{\Cal L}}$ simply as a line bundle,  in the case that ${\Cal L}$ is trivial as a line bundle on $X$.   Thus ${\Cal L}=\O_X$ with connection 
%$\nabla=d+\omega$,  where $\omega$ is an odd,  closed one-form on $X$.  Let us then denote ${\Cal L}$ by $\O_X^{\omega}$.  

Let $\Omega_{X,cl}^1$ denote the sheaf of closed one-forms.
Recall that the map  $\tilde d$ takes values in the sheaf  $\O_{\Delta}(1)$.   Furthermore, $\O_{\Delta}(1)=\Omega^1_X$.

\begin{Lem}\label{lem:closed}
The image of $\tilde d$ lies in $\Omega_{X,cl}^1$.
\end{Lem} 
(The next lemma implies that it is in fact all of $\Omega_{X,cl}^1$.)

\begin{pf}  Let $(z,\theta)$ be local coordinates on a neighborhood $\Cal U\subset X$.
Let $f\in\O_{\Delta}(\Cal U)$.
Then there is a one-form $\omega\in\Omega^1_X$ such that $f=\frac{\omega}{d\theta}$.  Write 
$d\omega=d\theta\wedge\alpha$,  $\alpha\in\Omega^1_X$.  Then $\tilde df=\alpha$.  Furthermore,
$0=d^2\omega=d\theta\wedge d\alpha$.  We can cancel $d\theta$, so $d\alpha=0$.

\end{pf}

%For a superscheme $Z$, let $\Cal N_Z\subset\O_Z$ denote the sheaf of nilpotents.

\begin{Lem}\label{lem:onto}
The map $\tilde d:\Cal N_{\Delta}\to\Omega_{X,cl}^1$ is surjective.\end{Lem}

\begin{pf}  Let $\alpha\in\Omega_{X,cl}^1$.  Then \dsp{\theta\frac{\alpha}{d\theta}\in\Cal N_{\Delta}} and \dsp{\tilde d (\theta\frac{\alpha}{d\theta})=\alpha}. 
\end{pf}

\begin{Rem}
Lemmas \ref{lem:closed} and \ref{lem:onto} imply the slightly weaker statement that there is an exact sequence
\begin{equation}\label{eq:weakexact}
0\to \O_{\hat X} \to \O_{\Delta} \overset{\tilde d}\to \Omega_{X,cl}^1\to 0 
\end{equation}
On the other hand,  it is known  \cite{BR1999} that the quotient of $\O_{\Delta}$ by $\O_{\hat X}$ is the Berezinian sheaf,  $Ber_{\O_{\hat X}}$.
We therefore have the corollary
\begin{Cor}\label{cor:ber}
$\Omega_{X,cl}^1\simeq Ber_{\O_{\hat X}}$ as $\O_{\hat X}$-modules.
\end{Cor}
\noindent This result seems to be new in this generality, although it is known for super Riemann surfaces where $X=\hat{X}$ \cite{GN1988, RSV1988}.

\end{Rem}

By lemmas \ref{lem:closed} and \ref{lem:onto} we have an exact sequence
\begin{equation}\label{eq:exact}
0\to (\Cal N_{\hat X})_0 \to (\Cal N_{\Delta})_0 \overset{\tilde d}\to (\Omega_{X,cl}^1)_1\to 0 
\end{equation}

\begin{Thm}\label{thm:class}
Let $\omega$ be an odd, closed one-form on $X$.  Let $c_\omega$ denote the image  of $\omega$ in
$H^1(\hat X,\Cal N_{\hat X})_0$ under the connecting homomorphism. Regarding $\widehat{\O_X^{\omega}}$ simply as a line bundle,  its class  in $H^1(\hat X,\O_{\hat X}^*)$ is $exp(c_{\omega})$.

\end{Thm}

\begin{pf} Pull $\omega$ back to $\Delta$,  giving the line bundle with connection  $\O_\Delta^{\pi^*\omega}$.   
Let $(z,\theta)$ be a local chart on a neighborhood $\Cal U\subset X$. 
Following the prescription in the proof of corollary \ref{cor:connecting},  decompose the constant function $1\in\O_{\Delta}$ as
$$1=\phi_0+\rho\phi_1$$
where $\hat\nabla_{\theta}(\phi_i)=0$.  
 Then $\widehat{\O_X^{\omega}}$ is trivialized on $\Cal U$ by 
 $\phi_0$.
 We have
 \dsp{\hat\partial_{\theta}=\frac 1{d\theta}\tilde d}.
 Then
 \begin{align}\label{eq:triv}
\phi_0&=1-\theta\hat\nabla_{\theta}(1)\\
&=1-\frac{\theta}{d\theta}\nabla_\theta(1)=1-\frac{\theta\omega}{d\theta}\\
&=exp(-\frac{\theta\omega}{d\theta})
\end{align}
Furthermore,
\begin{equation}\label{eq:}
\tilde d(\frac{\theta\omega}{d\theta})=\frac{d(\theta\omega)}{d\theta}=\omega
\end{equation}
\end{pf}

% There are  analogous results to theorem \ref{thm:class} in the setting of projected supercurves.  First a  lemma.

%\begin{Lem}\label{lem:pushforward}
%Let $\sigma:X\to X_0$ be a projected supercurve.   Then $\sigma_+$ takes lines bundles with connection to line bundles with connection. 
%\end{Lem}

%\begin{pf}
%Let ${\Cal L}$ be a line bundle with connection on $X$.  Let $\phi$ be a local trivialization in a neighborhood  on which there are local coordinates $(z,\theta)$.  Write $\phi=\phi_0+\theta \phi_1$,  where $\phi_i\in\sigma_+({\Cal L})$.   Write $\phi_i=f_i\phi$.  Then $f_0+\theta f_1=1$,  so $f_0$ is a unit.  Thus we may assume $\phi\in\sigma_+({\Cal L})$.  Then $\phi$ is a free generator of $\sigma_+({\Cal L})$.
%\end{pf}

Let us also note the following special case of corollary \ref{cor:connecting}.

\begin{Thm}\label{thm:pushforward}
Let  $\sigma:X\to X_0$ be a projected supercurve.  
Then

\smallskip

1.  The sequence
\begin{equation}\label{eq:pushexact}
0\to \Cal N_{X_0}\to\Cal N_X\overset d\to\Omega_{X/X_0,cl}^1\to 0
\end{equation}
is exact, where $d$ is the relative differential.\smallskip

2.  Let $\omega$ be an odd, closed one-form on $X$.  Let $\omega'\in \Omega_{X/X_0,cl}^1$ denote the image of $\omega$ under the natural map $\Omega_{X}^1\to\Omega_{X/X_0}^1$.
Let
$c_\omega$ denote the image  of $\omega'$ in
$H^1(X_0,\Cal N_{X_0})_0$ under the connecting homomorphism.   Then the class of $\sigma_+(\O_X^{\omega})$   in $H^1(X_0,\O_{ X_0}^*)$ is $exp(c_{\omega})$.
\end{Thm}

%\begin{pf}
%Working in local coordinates,  write a section of $\Omega_{X/X_0,cl}^1$ as $d\theta B$,  where $B$ is relatively closed.   Then we have $d\theta B=d(\theta B)$,  which establishes statement 1.

%
%For statement 2,  cover $X$ by neighborhoods $\Cal U_i\subset X$ with local coordinates  
%$(z_i,\theta_i)$.  
%Writing $\omega=dz_i A_i+d\theta_i B_i$.  Then the formula
%\begin{equation}\label{eq:decomp1}
%1=(1-\theta_i B_i)+\theta_i B_i
%\end{equation}
%gives the expansion of $1$ in powers of $\theta_i$ with coefficients in $\sigma_+(\O_X^{\omega})$.
%By the proof of lemma \ref{lem:pushforward},   $\sigma_+(\O_X^{\omega})$ is trivialized by $1-\theta B_i$.  Furthermore, applying the  differential relative ,  we have
%$d(\theta B_i)=d\theta B_i=\omega'$.  The result follows as in theorem \ref{thm:class}.
%\end{pf}

Our final result along these lines is a refinement of theorem \ref{thm:class} in the case that $\O_X^{\omega}$ is the pullback of  
the trivial bundle with connection on $X_0$.

\begin{Thm}\label{thm:class0}
Let  $\sigma:X\to X_0$ be a projected supercurve. Let $\Cal I\subset \O_{\hat X}$ denote the ideal sheaf of $X_0$ in $\hat X$ with respect to the corresponding imbedding, $\iota:X_0\to\hat X$.  Identify
$\Omega^1_{X_0}$ as a subsheaf of $\Omega^1_X$ by pullback.  Then
\smallskip

1.  The sequence
\begin{equation}\label{eq:pushexact}
0\to \Cal I \to\Cal I\Cal O_{X}\overset {\tilde d}\to\Omega_{X_0}^1\to 0
\end{equation}
is exact.\smallskip

2. Let $\omega$ be an odd one-form on $X_0$.  (Note that $\omega$ is necessarily closed for reasons of dimension.)  Let $c_\omega$ denote the image  of $\omega$ in
$H^1(\hat X,\Cal I)_0$ under the connecting homomorphism. Then the class of $\widehat{\O_X^{\omega}}$  in $H^1(\hat X,\O_{\hat X}^*)$ is $exp(c_{\omega})$.
\end{Thm}

\begin{pf}
In $(z,\theta)$ coordinates on $X$,  $\Cal I$ is generated by $\rho$.  Writing a section of
 $\Cal I\Cal O_{X}$ as $\rho\kern 1pt g(z,\theta)$,  we have $\tilde d(\rho g)=dz\frac{\partial g}{\partial \theta}$,  which shows that the image of $\Cal I\Cal O_{X}$ under $\tilde d$ is precisely $\Omega_{X_0}^1$.  The kernel is $(\Cal I\Cal O_{X})\cap\O_{\hat X}$,   which must be shown to coincide with $\Cal I$.
 One inclusion is obvious.  For the other inclusion,  take $\rho g\in \Cal I\Cal O_{X}$.  Then
 $\tilde d(\rho g)=0$ if and only if $\partial g/\partial\theta=0$,  which is to say $g\in\O_{X_0}$.  Then
 $\tau^{(z,\theta)}(g)=g+\rho\theta\frac{\partial g}{\partial z}$.  Thus $\rho g=\rho\tau^{(z,\theta)}(g)\in\Cal I$.  This completes the proof of statement 1.
 
 Statement 2 follows as in  theorem \ref{thm:class}.
\end{pf}

\subsection{Invertible sheaves in the complex topology}

%Up to this point,  there has been no need to specify whether we are working in the Zariski or complex topology.   

In
 this subsection, 
%  let us assume that $S=\Spec(\Lambda)$ where $\Lambda$ is supercommutative finite-dimensional $\BBB C$-algebra, and let us 
 we  work in the complex topology, where the Poincar\'e lemma is available.   Then the group of invertible sheaves on $X$ equipped with a flat connection is the group $H^1(X,\Lambda^*)$.  This group sees only the topology of $X$, and is therefore canonically isomorphic to  $H^1(\hat X,\Lambda^*)$.   In this way,  an
invertible sheaf with connection on $X$ induces an invertible sheaf with connection on $\hat X$,  by taking the same (constant) transition functions.

\begin{Prop}\label{prop:identity}
The identification of invertible sheaves with connection on $X$ and 
 invertible sheaves with connection on $\hat X$ given by ${\hat\pi}_+\pi^*$ coincides with the identity map on $H^1(X,\Lambda^*)$. 
 
 If $\sigma:X\to X_0$ is a projected supercurve,  then the identification of invertible sheaves with connection on $X$ and 
 invertible sheaves with connection on $X_0$ given by ${\sigma}_+$  also coincides with the identity map on $H^1(X,\Lambda^*)$. 
\end{Prop}
 
 \begin{pf} 
 We are given an invertible sheaf ${\Cal L}$ on $X$ with local trivializations $\phi_i$ on open sets $\Cal U_i$,  such that there are constants $c_{i,j}\in\Lambda^*$, such that $\phi_i=c_{i,j}\phi_j$.   The connection is then defined by the condition that $\nabla(\phi_i)=0$.   Letting $(z_i,\theta_i)$ be a coordinate system on $\Cal U_i$,  the local trivializations for $\hat{\Cal L}$ are
 $\tau^{(z_i,\theta_i)}(\phi_i)$.  By equation \eqref{eq:tau},
% Note that for all $\psi\in{\Cal L}$,
% $$\tau^{(z_i,\theta_i)}(\psi)=\frac{\nabla(\theta_i\psi)+dz\nabla_{\theta}(\psi)}{d\theta}$$
% In particular,  
 $\tau^{(z_i,\theta_i)}(\phi_i)=\phi_i$, so the transition functions are identical.
 
 The proof of the second assertion is the same.
 \end{pf}
 
 \section{Super elliptic curves}
 
 We illustrate the results of the previous section with a simple but nontrivial set of examples: super elliptic curves, i.e., supercurves of genus one \cite{R1995}.
 Let $\C^{1|1}$ be the trivial family of supercurves 
 $\Spec(\Lambda[z,\theta])$ over $S=\Spec(\Lambda)$.
 Let $X$ be the quotient of $\C^{1|1}$ by the discrete group $G \cong \Z \times \Z$ generated by the commuting morphisms
 $T(z,\theta) = (z+1, \theta)$ and $S(z,\theta)=(z + \tau + \theta \epsilon, \theta + \delta)$
 (not to be confused with the base scheme $S$).
 Here $\epsilon, \delta$ are odd elements of $\Lambda$ while $\tau$ is an even element satisfying
 ${Im}\, \tau_{{rd}} > 0$.
 This is a super elliptic curve whose underlying space $X_{{rd}}$ is the elliptic curve having  parameter $\tau_{{rd}}$.
 The corresponding $(1|0)$ curve $X_0$ is the quotient of $\C^{1|0}$ by the morphisms
 $T_0 (z) = z+1, \; S_0 (z) = z+\tau$.
 The dual curve $\hat{X}$ is easily computed as the quotient of $\Spec(\Lambda[u,\rho])$ by
 $\hat{T} (u,\rho) = (u+1, \rho)$ and 
 $\hat{S} (u, \rho) = (u + \tau + \epsilon \delta + \rho \delta, \rho + \epsilon)$,
 so that duality  exchanges $\epsilon$ with $\delta$ and changes $\tau$ to $\tau+\epsilon\delta$.
 $X$ is projected if $\epsilon=0$, injected if $\delta=0$, self-dual (a super Riemann surface) if
 $\epsilon=\delta$, and split if $\epsilon=\delta=0$.
 Only the self-dual case was considered in \cite{R1995}.
 The superdiagonal $\Delta$ is the quotient of $\C^{1|2}$ by
 $T_\Delta (z, \theta, \rho) = (z+1, \theta, \rho)$ and
 $S_\Delta (z, \theta, \rho) = (z+\tau+\theta \epsilon, \theta + \delta, \rho + \epsilon)$.
 
 We begin by determining the relevant cohomology of these curves.
 In case $\epsilon=\delta=0$, $H^0(X,\O_X)$ consists of functions $A(z)+\theta \alpha(z)$ 
 where $A$ and $\alpha$ are constants in $\Lambda$, since any nonconstant term in $A(z)$ or 
 $\alpha(z)$ of lowest degree in the generators of $\Lambda$ would give a nonconstant function on 
 $X_{{rd}}$, which is impossible.
 For general $\epsilon, \delta$, $H^0(X,\O_X)$ must be a submodule of this \cite{BR1999}.
 Clearly, these functions are invariant under the generator $S$ iff $\delta \alpha = 0$, so that
 $H^0(X,\O_X) = \Lambda | {ann} (\delta)$.
 $H^1(X,\O_X)$ is determined by Serre duality, since the dualizing Berezinian sheaf of $X$ is trivial, but a direct computation via group cohomology will provide more information, 
so we sketch it here \cite{R1995}.
 
 $H^1(X,\O_X) \cong H^1(G,\O) \cong H^1((S),\O^T),$ where $(S)$ is the cyclic subgroup generated by $S$, $\O$ are the functions on $\C^{1|1}$, and $\O^T$ are the $T$-invariant functions. 
 Geometrically this says that the cohomology of $X$ can be computed from the trivial cohomology of the cylinder arising from the quotient by $(T)$, by identifying its ends with $S$.
 A cocycle in $H^1((S),\O^T)$ assigns to the generator $S$ a $T$-invariant function 
 $A(z)+\theta \alpha(z)$, which is regarded as trivial if it can be written as
 $F(z,\theta) - F(S(z,\theta))$ for some $T$-invariant function $F(z,\theta) = f(z) + \theta \phi(z)$.
 This triviality condition implies
 \begin{align}  \label{trivial cocycle}
 A(z) &= f(z) - f(z+\tau) - \delta \phi(z+\tau),  \nonumber  \\
 \alpha(z) &= \phi(z) - \phi(z+\tau) - \epsilon f'(z+\tau) - \epsilon \delta \phi'(z+\tau).
 \end{align}
 Since all functions appearing are $T$-invariant, they have Fourier expansions of the form
 $A(z) = \sum_n A_n \exp 2 \pi i nz$, etc.
 Rewriting (\ref{trivial cocycle}) in terms of the Fourier components $A_n, \alpha_n$ shows that only the constant modes $A_0, \alpha_0$ can be nontrivial, in agreement with the expectation from Serre duality.
 For constant functions, (\ref{trivial cocycle}) immediately reduces to $A=-\delta \phi$.
 Thus we have $H^1(X,\O_X) = (\Lambda / \delta \Lambda) | \Lambda$.
 The cohomology of the dual curve $\hat{X}$ has the same form with $\epsilon$ replacing $\delta$, namely $H^1(\hat{X}, \O_{\hat{X}}) = (\Lambda / \epsilon \Lambda) | \Lambda$.
 
 With the usual exponential sheaf sequence, implying ${Pic}^0(X) = H^1(X,(\O_X)_0) / H^1(X,\Z)$, 
 this has the following interpretation.
 Line bundles of degree zero on $X$ can be specified by multipliers which are trivial for the cycle
 $T$ and of the form $\exp (A + \theta \alpha)$ for the cycle $S$, with $A$ and $\alpha$ even and odd elements of $\Lambda$ respectively.
Such a bundle is trivial when $\alpha=0$ and $A$ is a multiple of $\delta$.
For the dual curve, bundles having $\alpha=0$ and $A$ a multiple of $\epsilon$ are trivial.
Recall that proposition \ref{prop:identity} says that $\hat{\pi}_+ \pi^*$ relates bundles having the same constant transition functions on $X$ and $\hat{X}$.
This gives an example of a class in $H^1(X,\Lambda^*)$  defining the trivial bundle on $X$ and a nontrivial bundle on $\hat X$.  The existence of such classes was
pointed out in \cite{OR2002b}.

Our computation also allows us to determine which bundles in ${Pic}^0(X)$ admit flat connections, by
finding the image of $H^1(X,\Lambda)$ in $H^1(X,\O_X)$.
A cocycle for $H^1(G,\Lambda)$ assigns elements of $\Lambda$ to the generators of $G$, say
$T \mapsto -n, \;\; S \mapsto m$.
(The notation reflects the fact that this computation also determines the image of $H^1(X,\Z)$.)
To compare with our presentation of $H^1(G,\O)$, we subtract a trivial cocycle to set $n=0$, namely
$F(z,\theta)=nz$.
The result is $S \mapsto (m + n\tau) + \theta n \epsilon$.
That is, the bundles on $X$ admitting flat connections have $S$-multipliers 
$\exp (A + \theta \alpha)$ with
$\alpha$ a multiple of $\epsilon$.

We can similarly compute the cohomology of $\Delta$ in this example.
Global functions on $\Delta$ have the form $A + \theta \alpha + \rho \beta + \theta \rho B$ with
$A, B, \alpha, \beta \in \Lambda$.
We find that $H^0(\Delta, \O_\Delta)$ is the submodule of $\Lambda^2 | \Lambda^2$ given by the conditions $\delta B = \epsilon B = 0, \;\; \delta \alpha + \epsilon \beta = 0$.
Cocycles for $H^1(\Delta, \O_\Delta)$ have the same form, with the trivial ones generated by
$\delta \Lambda \cup \epsilon \Lambda$ and multiples of $\theta \epsilon - \rho \delta$.
Thus, for example, bundles on $X$ having multiplier $\exp A$ with $A \in \epsilon \Lambda$ would be trivial on $\hat{X}$ and lift to trivial bundles on $\Delta$; in addition bundles on $X$ having multiplier
$\exp \theta \epsilon$ and bundles on $\hat{X}$ having multiplier $\exp \rho \delta$ lift to the same
bundle on $\Delta$.

To illustrate theorem \ref{thm:class} we determine the closed one-forms on $X$; these have the form
$\omega = dz\,A + d \theta \, B$ where $A,B \in \Lambda$ and $G$-invariance requires $A \epsilon = 0$.
(The form $dz\, \delta + d \theta\, \theta \epsilon$ is also global, but not closed.)
Observe that $H^0(X,\Omega_{X,cl}^1) \simeq H^0(\hat{X},Ber_{\O_{\hat X}}) \simeq H^0(\hat{X},\O_{\hat{X}})$ as required by corollary \ref{cor:ber}.
Working through the proof to compute $\widehat{\O_X^\omega}$ we have
$\phi_0 = 1 - \theta B - \theta \rho A$.
Then the multiplier for $\widehat{\O_X^\omega}$ is $(\phi_0 \circ S) / \phi_0$, namely
$\exp -\delta (B + \rho A)$, which does indeed belong to the image of $H^1(\hat{X},\Lambda)$.
In the case $\delta=0$ when $X$ is injected, the transform of the trivial bundle is still trivial, but not in general.

To illustrate the other theorems, specialize to the case of $X$ projected, $\epsilon=0$.
Then theorem \ref{thm:class0} describes the transform of the pullback to $X$ of the trivial bundle with connection on $X_0$. 
Since the closed one-forms on $X_0$ have the form $\omega = dz\,A$, this is the special case $B=0$ of the result just obtained: the transform has multiplier $\exp -\delta \rho A$.

For theorem \ref{thm:pushforward}, begin with $\O_X^\omega$ where $\omega = dz\,A + d \theta \, B$ and there is no restriction on $A,B$ in this projected situation.
The image of $\omega$ in $\Omega_{X/X_0,cl}^1$ is $d \theta \,B$ and we have $\phi_0=1-\theta B$.
From the change in $\phi_0$ under $S$ we find the multiplier $\exp -\delta B$ for the direct image bundle on $X_0$.

% \subsection{Flat connections on the trivial bundle}
% Let $\omega$ be an odd, closed one-form on $X$.  Then the operator $d+\omega$ defines a flat connection on the trivial bundle $\O_X$.  Denote it by $\O^{\omega}$.    We then may transform $\O^{\omega}$ to obtain a line bundle with connection on $\hat X$.

%Some relationships between the Picard groups of $X$ and $\hat X$ were considered in \cite{BR1999}.
%In particular, a given set of constant ($\O_S$-valued) transition functions defines a line bundle on both $X$ and $\hat X$.
%Indeed, such bundles are related by our functor.

%\begin{Prop}
%Let $\F$ be an invertible sheaf on $X$, represented by a constant \v{C}ech cocycle in 
%$H^1(X,\O_X^\times)$, and viewed as a $\D_X$-module by taking $d$ as connection.
%Then $\hat\F$ is represented by the same cocycle in $H^1(\hat X, \O_{\hat X}^\times)$, with connection....
%\end{Prop}

%It was pointed out in \cite{OR2002b} that it is possible for $\hat \F$ to be trivial although $\F$ is nontrivial.
%The nontrivial connection $\hat \nabla$ in this case ``remembers" the nontriviality of $\F$.

\end{document}